\documentstyle{amsppt}
\magnification = 1200
\hcorrection{.25in}

\input psfig.tex

\document

\topmatter 

\title  {Grafting, harmonic maps, and projective structures on surfaces} 
\endtitle 

\rightheadtext {Grafting, harmonic maps, and projective structures}

\author Harumi Tanigawa  
\endauthor  

\address
Graduate School of Polymathematics,  Nagoya University, Nagoya 464-01, 
Japan 
\endaddress

\abstract 
Grafting is a surgery on Riemann surfaces introduced by Thurston
which connects 
hyperbolic geometry  and the theory of projective structures on surfaces. 
We will discuss the space of projective structures in terms of
the Thurston's geometric parametrization given by grafting. From
this approach we will prove that on any compact Riemann surface with genus
greater than
$1$ there exist infinitely many projective structures with Fuchsian
holonomy representations. In course of the proof it will turn out that
grafting is closely related to harmonic maps between surfaces.

\endabstract

\endsubjclass
\thanks {Research at MSRI is supported by NSF grant 
\#DMS--9022140}
\endthanks

\endtopmatter 
\head
1. Introduction
\endhead

A projective structure (or  a $\operatorname {\bold CP}^1$-structure)
on a surface is a coordinate system modelled on the
projective space $\operatorname {\bold CP}^1$ such that the transition maps
are projective homeomorphisms (hence the restriction of elements of
$\operatorname {PSL} (2,
\operatorname {\bold C})$). For an oriented closed surface $\Sigma_g$
of genus $g \ge 2$, it is well known that the space of projective
structures $P_g$ on $\Sigma_g$ is parametrized by the bundle of holomorphic
quadratic differentials on Riemann surfaces $\pi: Q_g \to T_g$ over the
Teichm\"uller space:  for each projective structure on $\Sigma _g$, take
the Schwarzian derivative of the developing map, and we have a quadratic
differential which is holomorphic with respect to the underlying complex
structure of the projective structure. As this parametrization is dealing
with projective or complex analytic mappings and manifolds, a lot of
researches have been developed from the viewpoint of complex analysis. (As
for this parametrization, see Hejhal \cite {H} for example.)\par
 The connection between projective structures on surfaces and hyperbolic
geometry was revealed by W. Thurston (unpublished). He showed that the
space $P_g$ of projective structures is parametrized by the product of
the Teichm\"uller space and the space of measured laminations. His idea is
to see a projective structure as a structure obtained  by bending a
hyperbolic 2-manifold in the hyperbolic 3-space along a measured geodesic
lamination.  Bending along a measured geodesic lamination is in some
sense conjugate to the earthquake deformation along the lamination (see
Epstein-Marden
\cite {EM} for detail). He also defined a surgery called grafting, which
is an equivalent concept with bending. 
\par
In this paper, we will study  projective structures and their
underlying complex structures from this geometric viewpoint. 
Especially, we will investigate the underlying complex structures of
projective structures with discrete holonomy representations
whose developing maps are not covering maps. The existence of such projective
structures was shown by Maskit \cite {Ma}, Hejhal \cite {He} and Goldman \cite
{G}, while it was unknown on which complex structure such projective structures
exist.
We will show
that on {\it any} complex structure on $\Sigma _g$ there are
infinitely many projective structures with Fuchsian holonomy
representations.  
To prove this fact, we will
define a mapping on the Teichm\"uller space to itself by grafting, which is 
conjugate  to the earthquake.  

We will prove our results in section 3 after describing
bending, grafting and the Thurston's parametrization theorem in section 2.
\par 
In course of  arguments, we will
see that harmonic maps are involved in grafting: when we consider a
projective structure as a bent hyperbolic structure, the bent surface
is a generalization of a pleated surface for the holonomy
representation, which is not necessarily discrete (see section 2). 
(In fact, when the holonomy representation is
discrete, the bent surface is a pleated surface of the quotient
3-manifold.)
On the other hand,
pleated surfaces in hyperbolic $3$-manifolds are the limits of the images of
harmonic maps by Minsky \cite {Mi2}.  We will see that the inverse of bending can be
seen as  mappings from  Riemann surfaces to the generalized pleated surfaces, so
that grafting is naturally related harmonic maps, in view of \cite
{Mi2} (see Remark 1  after Theorem 3.4).  

\par
   The author would like to thank Curt McMullen for his considerable help
and encouragement through this project. Most of this work was done at
Mathematical Sciences Research Institute, where the author
enjoyed various help by many people. Especially, she is very
grateful   to Michael Kapovich to whom she owes a lot on the proof
of the local injectivity of grafting, to William Thurston for his inspiring
explanation on the geometric parametrization of projective structures, 
and to Michael Wolf for useful and enjoyable discussions on
harmonic maps and the theory of measured laminations.
         
\head
2. Bending, grafting and geometric parametrization of projective structures
\endhead
In this section we sketch the Thurston's geometric
parametrization theorem. This geometric description of projective
structures is given by two equivalent concepts, {\it bending} or {\it
grafting}, which we will describe in this section.
Bending plays a role similar to that of pleated surfaces for hyperbolic
3-manifolds. Roughly speaking, bending is 
the way to see a projective
structure as a hyperbolic structure bent in the hyperbolic 3-space, and
grafting is the observation of bending on the sphere at infinity.
 
\subhead 
2.1. Thurston metric
\endsubhead
 We begin with a
metric introduced by Thurston which is a powerful tool to understand
projective structures.\par

Recall that every complex structure on a compact oriented surface
$\Sigma_g$ of genus $g$ admits a unique hyperbolic structure. This fact
provides two different approaches for Teichm\"uller theory: the
Teichm\"uller space $T_g$ is the space of complex structures and, at the
same time,  the space of hyperbolic structures on a compact surface
$\Sigma_g$. Now, for any complex structure $X \in T_g$ the
set of projective structures on $X$ are parametrized by the space of
holomorphic quadratic differentials on
$X$, which is a $3g-3$-dimensional complex vector space. As the complex structures
under these projective structures are all the same, the hyperbolic metric
does not distinguish them. 
The metric structure which characterizes a projective structure 
is defined by a very natural analogy of the definition of
hyperbolic metrics.
  
\definition {Definition 2.1 (Thurston (pseudo-)metric) } Let $M$ be a
$\operatorname{\bold CP}^1$-manifold. For each point $x\in M$ and
each tangent vector $v \in T_x M$, define the length of the vector $v$ by 
$$ t_M(v) = \inf_{f:\Delta \to M} \rho _\Delta (f^*v) $$
where the infimum is taken over all projective immersions $f: \Delta \to M$
with $f(\Delta) \ni x$ and
$\rho _\Delta$ is the hyperbolic metric on the unit disc $\Delta =
\{z\in \operatorname{\bold C} ;|z| < 1\}$. We will call the pseudometric
$t_M$ the Thurston pseudometric on $M$. If $t_M$ is non-degenerate it will
be called the Thurston  metric.
\enddefinition

Recall that the Kobayashi metric on a Riemann surface, 
which coincides with the hyperbolic metric if non-degenerate, is
defined by taking the infimum over all holomorphic immersions. 
 The
following properties are immediate consequences from the definitions of
Thurston  metric and Kobayashi hyperbolic metric.

\proclaim {Proposition 2.2} For a $\operatorname{\bold CP}^1$-manifold
$M$, let $k_M$ denote the Kobayashi pseudo-metric on $M$.
Then
\roster
\item $t_M \ge k_M$.
\item  If these metrics are non-degenerate on
$M$ and coincide at a non-zero tangent vector $v$ then these two metrics
coincide   on the entire tangent space $TM$.
\item For the  projective universal covering space $\tilde M$ of $M$,
$t_{\tilde M}$ descends to $t_M$ via the projective universal covering map
$\tilde M
\to M$.
\item If $t_M(v) \ne 0$ for a vector $v \in T_z M$ at a point $z \in M$
 then there is  a projective mapping $f : \Delta  \to M$ that attains
the minimum in the definition of $t_M(v)$. The mapping $f$ is
determined by $z$ uniquely up to precomposition of automorphisms of
$\Delta$.
\endroster
\endproclaim
In the following, we assume that the underlying complex structure of the
$\operatorname{\bold CP}^1$-manifold $M$ is hyperbolic, hence $t_M$
does not degenerate.\par 
 For convenience, we consider the Thurston metric on the universal
projective covering space $\tilde M$ rather than on $M$, as any
extremal mapping $f :\Delta \to \tilde M$ which realizes  the Thurston
metric at $z \in \tilde M$ is an embedding. \par
For each point $z \in \tilde M$ 
 the image
$f(\Delta)$ by an extremal mapping $f$ is a disc determined uniquely by $z$.
(Note that the terminology ``discs" makes sense in
$\operatorname{\bold CP}^1$-manifolds.) 
This disc is called the {\it maximal disc} for $z$. Let $D_z$ denote the
maximal disc for $z\in \tilde M$. Take a projective mapping $f$ 
on the upper half plane to $\tilde M$
realizing
the Thurston metric at $z$ and identify $D_z$ with the upper half plane
model of the hyperbolic 2-space $\operatorname{\bold H}^2$  via $f$. Then
we can compactify $D_z$ with the circle at infinity $\operatorname {\bold R}
\cup \{\infty\}$ of $\operatorname{\bold H}^2$.  
Let $\zeta \in \partial D_z$ be a boundary point. If the mapping $f :
\operatorname{\bold H}^2 \to \tilde M$ can be extended as a projective
map beyond $\zeta$, then $\zeta$ is identified with a point in the
frontier of $D_z$ in $\tilde M$. Otherwise, we call $\zeta$ an {\it
ideal boundary point}. Denote the set of all ideal boundary points of
$D_z$ by $\partial _\infty D_z$. Take the convex hull of $\partial
_\infty D_z$ with respect to the hyperbolic metric of $D_z ( =
\operatorname{\bold H}^2)$, and denote it by $C(\partial _\infty D_z)$.
It is easy to see that $\partial _\infty D_z$ consists of at least 2
points and there are three cases as follows (see Figure 1);
\roster
\item "{(i)}" $\partial _\infty D_z$ contains at least three points and $z$ is
      in the interior of 
      $C(\partial _\infty D_z)$.
\item "{(ii)}"  $\partial _\infty D_z$ contains at least three points and $z$ is
       in the frontier of $C(\partial _\infty D_z)$ in $D_z$.
\item "{(iii)}" $\partial _\infty D_z$ consists of two points and $z \in 
      C(\partial _\infty D_z)$.
\endroster
\midinsert
\hbox to\hsize{\hss\psfig{file=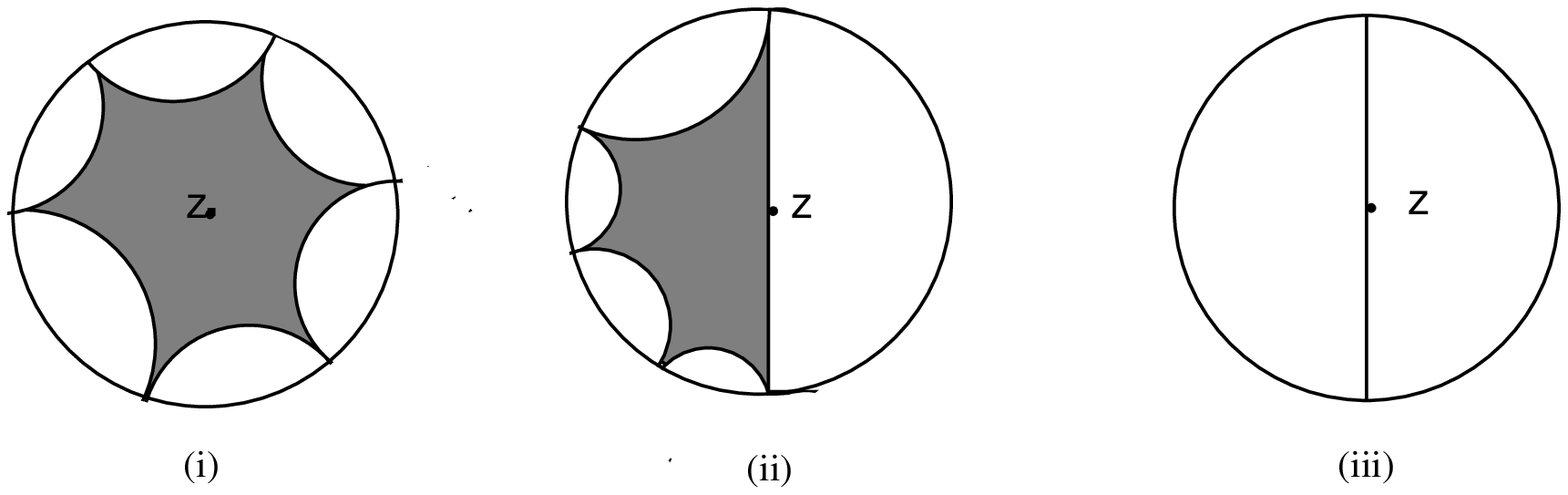}\hss}
\botcaption {Figure 1} The
point $z$ is in the convex hull of $C(\partial _\infty D_z)$.
\endcaption
\endinsert
We may assume that $0$ and $\infty$ are
ideal boundary points and
$z$ is on the imaginary axis.
In the first case, the Thurston metric
coincides with the hyperbolic metric
$|dz|/{\operatorname {Im} z}$ (on the upper half plane model of $D_z$)
near $z$. 
 
In the third case, Thurston metric at $z$ is equal to the flat
metric $|dz|/|z|$. In the second case, the hyperbolic metric and the
flat metric coincide  on the imaginary axis.\par
It is easy to see that 
$\tilde M$ is decomposed into 
the union of hyperbolic pieces and flat lines by the convex hulls of ideal
boundary points set
$C(\partial _\infty D_z)$ of maximal discs $D_z$.\par

\example {Example 2.3} Let $M = \tilde M$ be the union of two discs $D$
and $D'$ intersecting with angle $\theta \in [0, \pi)$ (Figure 2). For
convenience, we employ the model such that the intersecting two points
are $0$ and $\infty$. Let $S$ be the sector bounded by the ray
perpendicular to $\partial D$ and the ray perpendicular to $\partial D'$.
It is easy to see that for $z \in S$ the maximal disc for $z$ is the half
plane with boundary orthogonal to the ray through $z$ starting at $0$.
In this case the Thurston metric is equal to $|dz|/|z|$ on the ray. If
$z$ is outside of $S$ and contained in $D$ (resp. $D'$) then the maximal disc
for
$z$ is $D$ (resp. $D'$) and the Thurston metric near $z$ coincides with the
hyperbolic metric on $D$ (resp. $D'$). \par 
Therefore, Thurston
 metric is hyperbolic in $D-S$ and $D'-S$ and
flat in $S$. 
\endexample
\midinsert
\hbox to\hsize{\hss\psfig{file=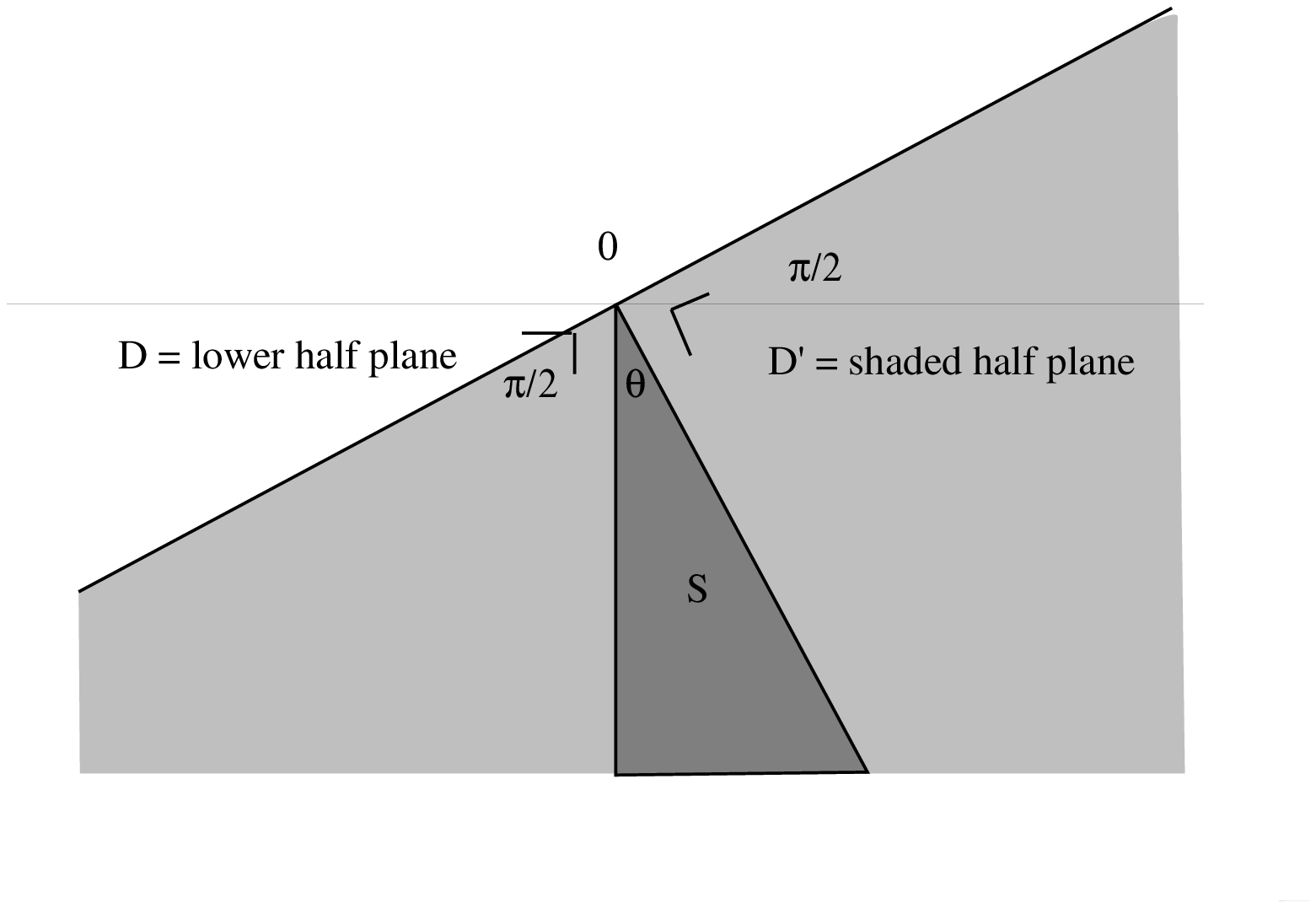}\hss}
\botcaption {Figure 2}
 A projective surface consists of hyperbolic pieces and flat pieces.
\endcaption
\endinsert
Note that in fact $\theta$ can be any positive number; if $\theta \ge
\pi$, we distinguish each sheet over the overwrapping region by
regarding the surface as $\{(re^{i\rho}, \rho) \in \operatorname {\bold C}
\times 
\operatorname {\bold R}:r \ne 0,  0 < \rho < \theta\}$.\par

\subhead
2.2. Bending a hyperbolic surface in $\operatorname {\bold H}^3$
\endsubhead

Next, we shall see that projective structures are obtained by bending
the hyperbolic 2-space $\operatorname {\bold H}^2$ in a locally convex way
in the hyperbolic 3-space $\operatorname {\bold H}^3$. In what follows,
we will denote by $CH(E)$ the convex hull of a subset $E$ in 
$\operatorname {\bold H} ^3 \cup \operatorname {\bold CP}^1$, where 
$\operatorname {\bold CP}^1$ is considered as a sphere at infinity of 
$\operatorname {\bold H} ^3$, to avoid mixing up the convex hull in
$\operatorname {\bold H} ^3$ with that in $\operatorname {\bold H} ^2$.
\par

We begin with a simple  example. We
will consider the Riemann sphere as the sphere at infinity of the
hyperbolic space $\operatorname{\bold H}^3$. Let
$D$ be a disc in the Riemann sphere. The convex hull $CH(D)$ of $D$ in   
$\operatorname{\bold H}^3$ is the half space bounded by the hyperbolic
plane $CH(\partial D)$. The nearest point projection $D \to CH(\partial
D)$ sends the hyperbolic structure of $D$ to the hyperbolic structure of
$CH(\partial D)$. On the other hand, the hyperbolic structure of $D$
coincides with the projective structure as a domain of $\operatorname
{\bold CP}^1$.  Hence in this case  the
projective structure on $D$ is given by the hyperbolic surface
$CH(\partial D)$ in $\operatorname{\bold H}^3$ with nearest point
projection. \par
Now, take a geodesic line $l \in CH(\partial D)$, fix an
orientation of $l$ and denote the left (resp. right) part of $CH(\partial D) - l$
by $\Delta^0$ (resp. $\Delta^1$). Take a positive number
$\theta$  (For simplicity, we temporarily assume that $\theta < \pi$) and 
rotate $\Delta^1$ along $l$  by angle $\theta$. Then we have a pleated
surface  $R$ as in Figure 3. 
 
\midinsert
\hbox to\hsize{\hskip.1in\hss\psfig{file=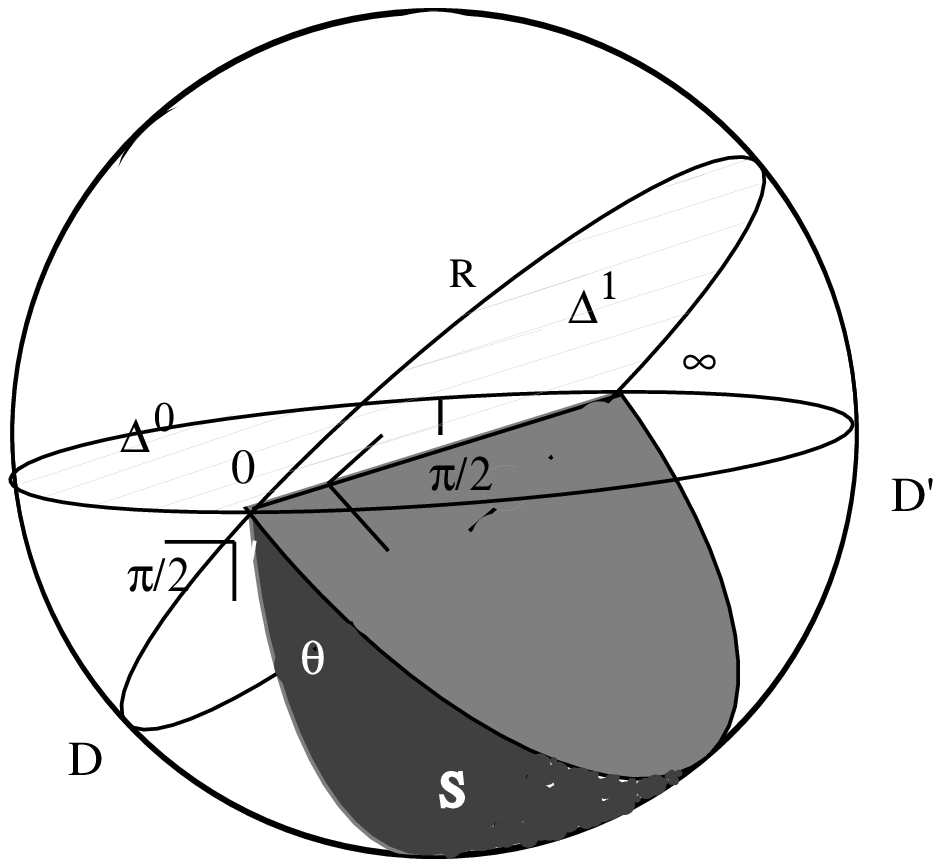}\hss\hskip-.1in}
\botcaption {Figure 3}
Bending a hyperbolic surface in the hyperbolic 3-space $\operatorname
 {\bold H}^3$ by angle $\theta$ produces a sector $S$ with angle
$\theta$ on the sphere at infinity, which is an Euclidean piece of the
projective surface
$\Omega  = D + D'$.
\endcaption
\endinsert

We will call this procedure  {\it bending the hyperbolic surface}
$CH(\partial D)$ {\it   along} $l$. \par
Now let us see what happens in the
sphere at infinity when we bend $CH(\partial D)$ along $l$. (Roughly
speaking, we get a new projective surface by
pushing the bent surface down to the sphere at infinity via the nearest
point projection.) \par

As we bend $CH(\partial D)$ in $\operatorname {\bold H}^3$ along  $l$,
$CH(\partial D)$ splits into two totally geodesic pieces, which are the
images of $\Delta^0$ and
$\Delta^1$. We denote the images by the same symbols
$\Delta^0$ and $\Delta^1$. For each of them,
there is a unique circle on the sphere at infinity whose convex hull in
$\operatorname {\bold H}^3$ contains the piece. For 
$\Delta^0$, the circle is the boundary of $D$. For $\Delta^1$, the circle
bounds the disc $D'$  intersecting with $D$ at the
endpoints of $l$ with angle $\theta$. 
Therefore, when we bend $CH(\partial D)$ in $\operatorname {\bold H}^3$
along  $l$ with angle $\theta$, the original projective surface
$D$ turns into
the  domain
$\Omega = D
\cup D'$. 
 This domain $\Omega$ has a projective structure as a domain of
the projective surface $\operatorname {\bold CP}^1$, which we observed
in Example 2.3. \par
We can reconstruct the pleated surface $R$ from $\Omega$ in the following
way. Remember that we saw in section 2.1 that for each $z \in \Omega$
there is a unique maximal disc $D_z$. For each $z \in \Omega$, take the 
convex hull of the circle $\partial D_z$ in $\operatorname {\bold H}^3$.
Then send each point in the convex hull $C (\partial _\infty D_z)$ (defined in
section 2.1) of
$\partial _\infty D_z$ in the hyperbolic surface $D_z$ by the nearest
point projection to the convex hull of $\partial _\infty D_z$ in 
$\operatorname {\bold H}^3$.  
Recall that we saw in Example 2.3  that
$\Omega$ is decomposed into hyperbolic pieces $D-S$ and $D'-S$ and a flat
piece $S$ with respect to Thurston metric $t_\Omega$. 
Then by the nearest point projection, $D-S$ (resp. $D'-S$) is mapped  to
$\Delta^0$  (resp. $\Delta^1$) isometrically. As for the sector $S$, each flat line
connecting $0$ and $\infty$ is mapped to $l$ isometrically.
Thus the image of $\Omega$ is the pleated surface $R$, and the above
mapping $\Omega \to R$ is the inverse of the procedure getting the
projective structure $\Omega$ from the pleated surface $R$.\par

Thus the procedure bending  $CH(\partial D)$ in $\operatorname {\bold H}^3$
along a geodesic is equivalent to
`grafting'
a flat part $S$ into  the hyperbolic structure on $D$. 
\par
As before, note that we do not have to restrict $\theta$ to be smaller
than $\pi$: if $\theta \ge \pi$, distinguish overwrapping sheets.
\par
Now we proceed to the case with a group action. Let $\Gamma$ be a
co-compact Fuchsian group acting on $\operatorname{\bold H}^2$. Embed
$\operatorname{\bold H}^2$ in $\operatorname{\bold H}^3$ as a totally
geodesic surface. Let $X$ denote the hyperbolic surface
$\operatorname{\bold H}^2/\Gamma$. Take a simple closed geodesic curve
$\gamma$ on $X$. The lift of $\gamma$ on $\operatorname{\bold H}^2$ is a 
$\Gamma$-invariant set of geodesic lines. We can bend
$\operatorname{\bold H}^2$ along each of these geodesics with angle
$\theta$ step by step (see Epstein-Marden \cite {EM}).
In each step, on the sphere at infinity, we have a new projective
surface with a  grafted part to the preceding step, as we did in the
preceding example.   (In each step distinguish the overwrapping sheets, if
any, as we did in Example 2.3.)
Then we end up with a simply
connected projective manifold 
$\tilde M$ spread over the sphere at infinity, which is partly hyperbolic
and partly flat. \par

In view of the construction of $\tilde M$, it is
easy to see that there is a projective automorphism group $\tilde \Gamma$ acting on
$\tilde M$ isomorphic to $\Gamma$, hence in particular, $\tilde M/\tilde \Gamma$ is
homeomorphic to
$X$.
To consider $\tilde M$ as spread over the Riemann sphere as above is to map $\tilde
M$ to
$\operatorname {\bold CP}^1$ via the developing map and we have the holonomy
representation
$\chi :\Gamma \to
\operatorname {PSL(2, \bold C)}$. Then the above bending procedure is 
given by an equivariant map 
 from $\operatorname {\bold H}^2$ to $\operatorname {\bold H}^3$
with respect to $\Gamma$ and the holonomy representation which is
bent along the bending locus and
isometric elsewhere. 
 
Indeed, it is known that we can write down
the holonomy representation $\chi
:\Gamma \to
\operatorname {PSL(2, \bold C)}$ in terms of bending.
(We omit   the formulae. See
\cite {EM Chapter 3} for details. There, the homomorphism is
called the quakebend homomorphism.)
\par

It is also known that when the weighted simple closed curves converge to a
measured lamination, the equivariant maps converge  to the equivariant map
bent along the measured lamination and it defines the corresponding
projective structure.

  \par

See Epstein-Marden [EM] for detail.
 \par
\subhead
2.3. Grafting along a simple closed curve
\endsubhead
Grafting is the way to see the above
procedure directly on the quotient surfaces $X = \operatorname {H}^2 /
\Gamma $ and $M = \tilde M / \tilde \Gamma$ as in the following way. 
\par
We provide two types of $\operatorname {\bold CP}^1$-manifolds which
we will paste together. Let $X$ and $\gamma$ be as in section 2.2.
First, take
the lower half plane model of
$\operatorname{\bold H}^2$ such that the geodesic line $\{iy; y < 0\}$
is  one of the component of the lift of
$\gamma$. Let $g(z) = e^{l(\gamma)}z $ be the generator of the
stabilizer of   
$\{iy; y < 0\}$ in $\Gamma$, where $l(\gamma)$ denotes the hyperbolic
length of
$\gamma$ on $X$. 
Next, take the sector $\{z=re^{i\rho}; 0<r<\infty, 0 \le
\rho \le \theta\}$ equipped with the projective structure as a domain of
$\operatorname{\bold CP}^1$. The group $<g>$ generated by
$g$ acts on this sector as a projective automorphism. Taking the
quotient we get a flat annulus $A_\theta$ with height $\theta$ and
contour $l(\gamma)$. \par
Now we cut $X$ along $\gamma$ and paste each side of
$\gamma$ to one of the boundary component of $A_\theta$ as Figure 4, in
such a way as the length parameters of pasting sides match and that the pair
of points which is identified in $X$ are connected by segments in
$A_\theta$ orthogonal to the boundary. \par
\midinsert
\hbox to\hsize{\hss\psfig{file=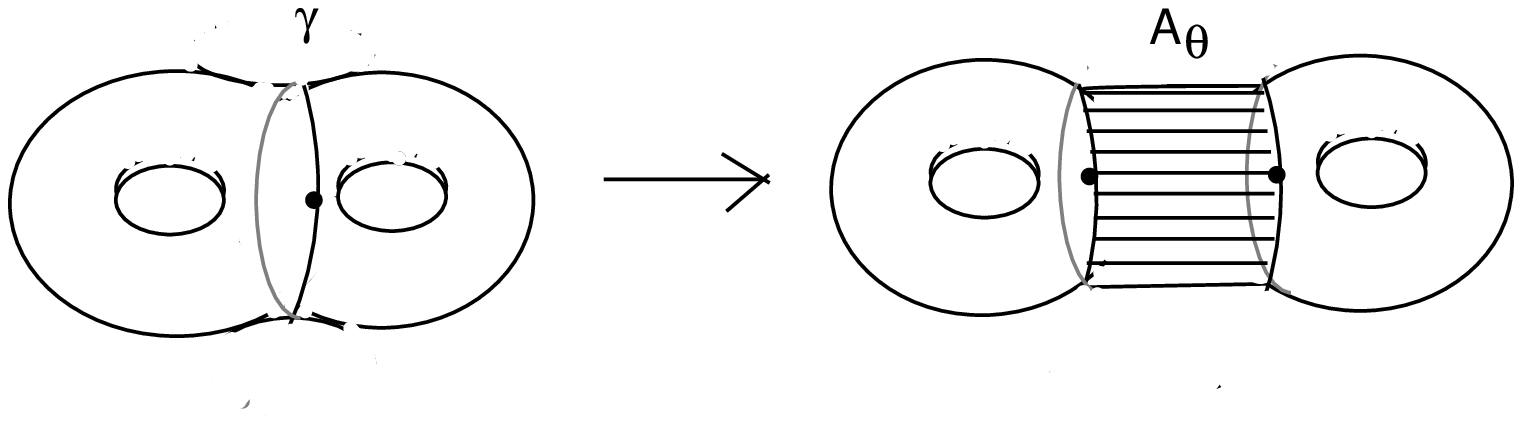}\hss}
\botcaption {Figure 4}
 Grafting a flat annulus of height $\theta$ to $X$ along $\gamma$.
\endcaption
\endinsert
Note that the hyperbolic structure of $X$ and the
projective structure of $A_\theta$ match on the pasting locus. Therefore,
the above pasting process yields a new $\operatorname {\bold
CP}^1$-structure preserving the original projective structures of $X$ and
$A_\theta$. We call this surgery {\it grafting a flat annulus of height
$\theta$ to $X$ along
$\gamma$}, or {\it grafting $\theta \gamma$ to $X$}, and denote the
resulted $\operatorname {\bold CP}^1$-structure by 
$Gr_{\theta
\gamma} (X)$. \par
 Note that the metrics on the hyperbolic surface $X$ and the flat annulus
$A_\theta$ also match on the pasting locus and the resulted surface
$Gr_{\theta
\gamma}(X)$ is equipped with a metric which is partly hyperbolic and partly
flat. It is easy to see that this metric is equal to the Thurston metric on
$Gr_{\theta 
\gamma}(X)$.  It follows that the underlying complex structure of
$Gr_{\theta  \gamma}(X)$ differs from $X$ by Proposition 2.2 (4), unless
$\theta = 0$.
\par
It
is also easy to see that this projective structure 
$Gr_{\theta \gamma}(X)$ has a projective universal covering space
$\tilde M$, which we obtained by bending in section 2.2.  
\subhead
2.4. grafting a general measured lamination and the parametrization theorem
\endsubhead
 
Let $P_g$ denote the set of all projective structures on the oriented
closed surface $\Sigma_g$ of genus $g$. Then as we have seen above, the
grafting operation gives a mapping
$$Gr : T_g \times \operatorname {\bold R}_+ \times \Cal {S} \to P_g,$$
which sends each $(X, \theta, \gamma) \in T_g \times \operatorname {\bold
R}_+
\times \Cal {S}$ to the projective structure obtained by grafting a flat
annulus of height $\theta$ along the hyperbolic geodesic in the homotopy
class of $\gamma$ to the hyperbolic surface $X$, where $\Cal {S}$ denotes
the set of homotopy classes of simple closed curves.  Now we can state 
Thurston's parametrization theorem.

\proclaim {Theorem 2.4 (Thurston)} The map $Gr$ extends to a
homeomorphism of
$T_g \times \Cal{ML}$ onto $P_g$, where $\Cal {ML}$ denotes the space of
measured laminations on $\Sigma_g$.
\endproclaim
\demo {Sketch of the proof} We have already seen that for any measured
lamination $\mu = \theta \gamma$  supported on a simple closed curve
$\gamma$ and for any hyperbolic structure $X \in T_g$ grafting an
annulus with height $\theta$ yields a projective structure. 
The mapping $Gr : T_g \times \operatorname {\bold R}_+ \times \Cal {S} \to
P_g,$ is continuously extended to $T_g \times \Cal {ML}$, as bending is defined for
any measured lamination and depends on the lamination continuously.
See \cite {EM} for detail. \par
We shall describe the inverse correspondence:
 $P_g \to T_g \times \Cal {ML}$. By the  arguments in the preceding
sections, it suffices to show that any projective structure on $\Sigma_g$ is
obtained from the  bending procedure defined with an equivariant map
$\operatorname {\bold H} ^2 \to \operatorname {\bold H} ^3 $, bent along a
measured lamination and isometric elsewhere. \par
 Given a projective
structure on
$\Sigma_g$ take its projective universal covering $\tilde M$ and fix its
developing map. Begin with an open set  $U$ in $\tilde M$ small
enough so that the developing map restricted to $U$ is homeomorphic.
For each point $z$ in $U$, take the maximal disc $D_z$ for $z$ and embed it
into the Riemann sphere via the developing map. We identify $D_z$ with its
image. Then take the convex hull
of the circle $\partial D$ in $\operatorname {\bold H} ^3$ and denote
it by $R_z$. $R_z$ is a totally geodesic disc isometric to $D_z$ with
respect to the hyperbolic metrics on them via the nearest point projection.
  Now, as in section 2.1, take the convex
hull
$C(\partial _\infty D_z)$ of the ideal boundary points of $D_z$. Then
as we did in section 2.2 for the simple case without group action, send
$C(\partial _\infty D_z)$ into $R_z$
via the nearest point projection between  $R_z$ and
$D_z$. Denote the image of $C(\partial _\infty D_z)$
by $P_z$. If 
$C(\partial _\infty D_z)$ is of the type (i) in section 2.1, then $P_z$ is a
convex domain of $R_z$ which is the convex hull of $C(\partial _\infty
D_z)$ in $\operatorname {\bold H}^3$. If
$C(\partial _\infty D_z)$ is of the type (iii), then $P_z$ is the hyperbolic
line in $R_z$ connecting the two points in $C(\partial _\infty D_z)$. In
any case,
$P_z = P_{z'}$ for every $z' \in C(\partial _\infty D_z)$. Now, 
$\cup _{z \in U} P_z$ 
is a piece of a pleated surface in $\operatorname {\bold H} ^3$:
there are a subset $V$ of $\operatorname {\bold H}^2$ and
a mapping from $V$  to $\cup _{z \in U} P_z \subset \operatorname {\bold H} ^3$,
such that for each point $w \in V$
there is a straight line in $V$ which is mapped isometrically to a
hyperbolic line in $\operatorname {\bold H} ^3$.
This piece of pleated surface defines locally the bending which gives the  
projective structure of $\cup _{z \in U} D_z$. 
(Here, $\cup _{z \in U} D_z$ is equipped with the projective structure as
a domain of $\tilde M$.) 
  
Beginning with $U$ and continuing this procedure, 
  it is easy to get an equivariant mapping $\operatorname {\bold H}^2$
to $\operatorname {\bold H}^3$, defining the bending which produce the
projective surface $\tilde M$.
 See
\cite {KT} for detail.

\enddemo
  
\head
3. Grafted structures on surfaces
\endhead

Now, we are ready to discuss projective structures in terms of the
geometric parametrization. Given a measured lamination $\mu$, let $gr_\mu
(X)$ stand for the underlying complex structure of the projective
structure
$Gr_\mu (X)$ for $X \in T_g$.  For any fixed $\mu$, this
assignment gives a mapping 
$gr_\mu : T_g  \to T_g$. We shall call this mapping {\it the grafting map
defined by} $\mu$.
  
\par
First, we recall some facts about projective structures with Fuchsian
holonomy representations. On any complex structure $X \in T_g$ there is
a unique projective structure whose projective universal covering space 
is projectively equivalent to the hyperbolic 2-space $\operatorname {\bold
H}^2$, namely, the hyperbolic structure. The holonomy representation 
of this projective structure is a
Fuchsian group $\Gamma$ acting on $\operatorname {\bold H}^2$ with
quotient manifold $X = \operatorname {\bold H}^2 / \Gamma$. An `exotic'
projective structure with Fuchsian holonomy representation whose
developing map is not a covering map was first constructed by Maskit
\cite {Ma}. Hejhal \cite {He} and Goldman \cite {G} made more topological
and geometric approach to such projective structures. The following
characterization of projective structures with Fuchsian holonomy
representations was given by Goldman.
\proclaim {Theorem 3.1 (Goldman \cite {G})}
A projective structure given by $(X, \mu) \in T_g \times \Cal {ML}$
has a Fuchsian holonomy representation if and only if $\mu$ is an
integral point of $\Cal {ML}$. Here, a measured lamination $\mu$ is
called an integral point if it is of the form $\mu = \sum 2 \pi m_i
\gamma_i $ with a disjoint union of nontrivial simple closed geodesics
$\{\gamma_i\}$ and a set of positive integers $\{m_i\}$.
\endproclaim

Note that given a projective structure determined by a pair 
$(X, \mu) \in T_g \times \Cal {ML}$, 
the underlying complex structure  of $Gr_\mu (X)$ is hardly expressed by 
$X$ and $\mu$, unless
$\mu = 0$. So far, in particular, it is unclear on which complex structures
there exist projective structures with Fuchsian holonomy representations
other than the hyperbolic structures. Our main result shows that on {\it
any} complex structure and any integral point $\mu \in \Cal {ML}$ there is a
unique projective structure with Fuchsian holonomy representation which
is obtained by grafting
$\mu$ to some hyperbolic structure $X \in T_g$:

\proclaim {Theorem 3.2} For any integral point $\mu \in \Cal {ML}$, the
grafting map $gr_{\mu} : T_g \to T_g $ is a real analytic homeomorphism.
\endproclaim

Before proving this theorem, let us interpret it in terms of
the parametrization  of $P_g$  (the space of projective structures) by the
bundle of holomorphic quadratic differentials on Riemann surfaces
$\pi : Q_g \to T_g$. This parametrization is given in the following way:
for each projective structure, take the Schwarzian derivative of the
developing map, where the Schwarzian derivative of a locally univalent
meromorphic function $f$ is defined by $(f''/f')' - 1/2 (f''/f')^2$. Then
the Schwarzian derivative is a quadratic differential on the surface which
is holomorphic with respect to the complex structure under the projective
structure.  
(See Hejhal \cite {He}, for example, for detail.)
 The canonical projection $\pi : Q_g \to T_g$ sends each
projective structure to its underlying complex structure. Let
$K
\subset Q_g$ is the set of projective structures with discrete holonomy
representations. For $X \in T_g$, let $Q(X)$ and $K(X)$ denote the
fibers over $X$ of $\pi : Q_g \to T_g$ and $\pi |K : K \to T_g$
respectively. For every $X \in T_g$, the interior point set
$\operatorname {int} K(X)$ has a component containing $0$, which
coincides with the Bers slice. Theorem 3.2 implies the existence of
components of $\operatorname {int} K(X)$  other than the Bers slice: 
\proclaim {Corollary 3.3} On any complex structure $X \in T_g$, there
are infinitely many components of $\operatorname {int} K(X)$.
\endproclaim

\demo {Proof of Corollary 3.3}
Fix an integral point  $\mu \in \Cal {ML}$  and  a hyperbolic
structure $X'$. The projective structure $Gr_\mu
(X')$ has a Fuchsian holonomy group $\Gamma'$.
For any Beltrami differential $\tau$ for $\Gamma'$ on the Riemann
sphere $\hat {\operatorname {\bold C}}$, we can take a quasiconformal
deformation of the projective structure $Gr_\mu (X')$ by $\tau$ (cf.
\cite {ST}): let
$f^\tau$ denote the quasiconformal homeomorphism of 
$\hat {\operatorname {\bold C}}$ with Beltrami differential $\tau$ fixing
$0,1$ and $\infty$.  Then $\gamma ^\tau = f^\tau \circ \gamma \circ
(f^\tau)^{-1}$ is a M\"obius transformation for every
$\gamma
\in
\Gamma'$ and
$\Gamma ^\tau = f^\tau
\Gamma' (f^\tau)^{-1}$ is a quasifuchsian group. As $\gamma ^\tau \circ
f^\tau = f^\tau \circ \gamma$, we have another projective structure by
replacing the local coordinate system $\{ (U, \phi) \}$ of $Gr_\mu (X')$
to $\{ (U, f^\tau \circ \phi) \}$ with holonomy representation
$\Gamma^\tau$. It is easy to see that this new projective structure
depends only on the equivalence class of $\tau$ (cf. \cite {ST}).
Therefore, we have an open set $QF(\mu)$ of $K$ consisting of all
quasiconformal deformations of $Gr_\mu (X')$. Note that if $\Gamma
^\tau$ is a Fuchsian group, the projective structure defined by the
quasiconformal deformation of $Gr_\mu (X')$ with $\tau$ is equal to the
projective structure $Gr_\mu (X^\tau)$, where $X^\tau$ is the hyperbolic
surface obtained by the quasiconformal deformation of $X'$ with the
Beltrami differential $\tau$. Hence by Theorem 3.2 the restriction $\pi |
QF(\mu) : QF(\mu) \to T_g$ is surjective for each integral point $\mu \in
\Cal {ML}$. Therefore, the corollary follows if we show $QF(\mu) \cap
QF(\nu) = \emptyset$ for any two different integral points $\mu$ and
$\nu$.  To see this, take the inverse image of the limit set $\operatorname
{\bold R}\cup \infty$ of $\Gamma'$ via the developing map on the universal
cover of the $\operatorname {\bold CP}^1$-manifold $Gr_\mu (X')$. Then 
the inverse image
 descends to a
disjoint union of curves on $Gr_\mu (X')$. If the integral point
$\mu$ is of the form $\mu = \sum 2n_i\pi \gamma_i$ for  integers $\{n_i\}$
and  simple closed curves $\{\gamma_i\}$, 
then the inverse image of the limit set of $\Gamma'$ descends to the union
of 
$2n_i$ curves each of which is homotopic to 
$\gamma_i$. On the other hand, it is easy to see that any quasiconformal
deformation of the projective structure $Gr_\mu (X')$ maps this
system of curves quasiconformally (the image depends only on the
equivalence class of the Beltrami differential). Therefore, the
homotopy class of these system of curves characterizes the open set
$QF(\mu)$. Hence $QF(\mu) \cap QF(\nu) = \emptyset$ for
any two different integral points
$\mu$ and $\nu$. \qed
\enddemo

\par
\remark {Remark} It was shown by Maskit \cite {Ma} that there exists
{\it some}
$X$ such that $\operatorname{int} K(X)$ (the interior of $K(X)$ in $Q(X)$)
has {\it
some}   components other than the Bers slice. 
 In \cite {ST} we discussed on $\operatorname {int
} K (X)$ for such $X$ (i.e. assuming the existence of such components on
$X$), where we showed that any component of
$\operatorname {int} K(X)$  is  a component of $QF(\mu) \cap Q(X)$ for an integral
point
$\mu \in \Cal {ML}$. What we have shown in the above corollary is that 
$QF(\mu) \cap Q(X)$ is a non-empty open set for every
complex structure $X$ and every integral point $\mu$.
\endremark 
\bigskip
{\it Proof of Theorem 3.2.}
To prove Theorem 3.2, it suffices to show for an integral point $\mu$
\roster
\item $gr_\mu : T_g \to T_g$ is a proper mapping,
\item $gr_\mu : T_g \to T_g$ is a local diffeomorphism, and
\item $gr_\mu : T_g \to T_g$ is real analytic.
\endroster

\demo {Proof of (1)}The following theorem enables us to show that  for {\it
any} measured lamination
$\mu$ (not necessarily an integral point), the grafting map $gr_\mu : T_g
\to T_g$ is a proper map. 

\enddemo
\proclaim {Theorem 3.4}
Let $X$ be a hyperbolic surface and $\mu$ be a measured lamination. Let
$h : gr_\mu(X) \to X$ denote the harmonic map with respect to the
hyperbolic metric on $X$ and $\Cal {E} (h)$ be its energy.
(Remember that the harmonic map between surfaces
depends on the metric on the target surface but only on the conformal
structure on the source surface.) 
Then
$$\frac {1}{2} l_X(\mu) \le \frac {1}{2}\frac{l_X(\mu)^2}{E_{gr_\mu
(X)}(\mu)} \le 
\Cal {E} (h) \le \frac {1}{2} l_X(\mu) + 4\pi (g-1),$$
where $l_X(\mu)$ is the hyperbolic length of
$\mu$ on $X$ and $E_{gr_\mu (X)}(\mu)$ is the extremal
length of $\mu$ on the grafted surface
$gr_\mu(X)$.
\endproclaim
\demo {Proof of Theorem 3.4} For simplicity, we abbreviate $Y = gr_\mu
(X)$.  First, assume that
$\mu$ is supported on a simple closed curve, $\mu = \theta \gamma$. Then
the projective structure $Gr_\mu (X)$ consists of hyperbolic piece(s) whose
union is identified with $X$ and a flat annulus $A_\theta$. We will use
this geometric structure on $Y$. We define a mapping
 $f: Y \to X$ by collapsing the annulus $A_\theta$ to the geodesic curve
$\gamma$ on the hyperbolic surface $X$ along the flat structure (i.e.
translating each point of $A_\theta$ to $\gamma$ along the segment
perpendicular to
$\gamma$) and
 sending the hyperbolic pieces of $Y$ isometrically on the corresponding
domains on $X$. Then $f$ is among the competitive mappings for the
harmonic map $h: Y \to X$. As $f$ is isometric on the hyperbolic pieces,
the contribution of this part for the total energy is the hyperbolic
area of $X$. On the flat annulus $A_\theta$, the direction parallel to
the geodesic
$\gamma$ and the direction of the segment orthogonal to $\gamma$ forms
an orthogonal frame in $A_\theta$. The length of the
former direction is preserved by $f$, while the image of the latter
direction degenerates. Therefore, the contribution to the total energy
of the flat part is $ 1/2 \theta l_X(\gamma) = 1/2 l_X (\mu)$. Hence
we have
$$ \Cal {E} (h) \le \Cal {E} (f) \le \frac {1}{2}  l_X (\mu) + 4g-4.$$
On the other hand, by  the left part of Minsky's inequality \cite {Mi1,
Theorem 7.2},
$$\frac {1}{2}\frac{l_X(\mu)^2}{E_Y
(\mu)}\le \Cal {E}(h).$$
Note that the extremal length of $\gamma$ in $Y$ is not greater than that
in
$A_\theta$, that is, $ l_X(\gamma)/\theta$. It follows that 
$$ \frac {1}{2} l_X(\mu) = \frac {1}{2} \theta l_X(\gamma) = \frac
{1}{2}\frac {l_X(\gamma)^2} { l_X(\gamma)/\theta} 
\le
\frac {1}{2}\frac{l_X(\gamma)^2}{E_Y (\gamma)}
= \frac {1}{2}\frac{l_X(\mu)^2}{E_Y (\mu)}.$$
We have shown the inequality in the statement of Theorem 3.4 in the case
$\mu$ is supported on a simple closed curve.
For a general measured lamination, we approximate $\mu$ by a sequence of
measured laminations each of which is supported on a simple closed curve.
The inequality follows from the continuity of the hyperbolic length of
measured laminations on $X$ and the continuity of grafting with respect
to measured laminations.
\qed
\enddemo

\bigskip

Now we prove the properness of $gr_\mu: T_g \to T_g$ from Theorem 3.4.
When a sequence of points in $T_g$ leaves any compact set eventually,
we will say `the sequence tends to infinity' for simplicity.
We have to show for any sequence $\{X_n\}$ tending to infinity the image
$\{gr_\mu(X_n)\}$ also tends to infinity. Denote
$Y_n = gr_\mu(X_n)$ for simplicity.
By taking a subsequence if necessary, we may assume that either
\roster
\item "{(i)}" $\sup_n l_{X_n}(\mu) < \infty$, or
\item "{(ii)}" $\lim _{n \to \infty} l_{X_n}(\mu) = \infty.$
\endroster

In the case (i) we show that $\{Y_n\}$ tends to infinity by contradiction.
Assume that $\{Y_n\}$ stays in a compact set of $T_g$. As
$X_n$ tends to infinity, the energy of the harmonic map  $h_n : Y_n \to
X_n$ tends to infinity by a result of M. Wolf [W, Proposition 3.3]. This
contradicts the assumption that $l_{X_n} (\mu)$ is uniformly bounded,
considering the right inequality of Theorem 3.4.\par
In the case (ii), by Theorem 3.4,
$$\lim _{n\to \infty} E_{Y_n}(\mu) = \lim _{n \to \infty} (l_{X_n} (\mu ) +
O(1)) = \infty.$$ Therefore, $Y_n$ tends to infinity.
\qed
\bigskip
\remark {Remark 1 (Collapsing the grafted part is close to the harmonic
map)} In the above proof of  Theorem 3.4, we showed that 
the difference between the total energy of
the annulus collapsing map $f: Y \to X$  and that of the harmonic map
$h : Y \to X$ is bounded by a universal constant depending only on the
genus $g$.  Therefore, we can say that
$f$ is close to the harmonic map when
$l_X(\mu)$ is large,
as the harmonic map between a pair of hyperbolic
surfaces is unique by a result of Hartman \cite {Ha}.
Here we exhibit an 
intuitive explanation for this phenomenon 
\par  
First, note that the grafted part occupies a large portion on the
entire surface  when
$l_X(\mu)$ is very large, in view of the Thurston metric on
$Gr_\mu (X)$. To collapse this large part likely results in
`significant stretch in this direction'. 
In general, the direction of `maximal stretch' of any kinds of
extremal mappings (e.g. Teichm\"uller mappings, extremal Lipschitz maps, or
 harmonic maps)
plays the key role to measure the difference between two surfaces: 
Kerckhoff \cite {Ke} showed that the Teichm\"uller distance between two
Riemann surfaces is described by the ratio of the extremal lengths of
the direction of maximal stretch of the Teichm\"uller mapping.
Similar results for Lipschitz maps are proved by Thurston
\cite {Th2}, and for harmonic maps, by Minsky
\cite {Mi1} and \cite {Mi2}.  Now, as for grafting, it is natural to pay
attention to harmonic maps to compare the grafted surface with the original
surface for the following reason. \par
Recall (see section 2.2) that grafting a measured lamination $\mu$ to a
hyperbolic surface
$X$ is equivalent to bending which is realized by the equivariant map $g
:
\operatorname {\bold H}^2
\to \operatorname {\bold H}^3$, with respect to the Fuchsian group $\Gamma$
with $\operatorname {\bold H}^2 / \Gamma = X$ and the holonomy
representation of $Gr_\mu (X)$, which is bent along the lift of $\mu$ and
isometric elsewhere. This is a generalization of a pleated surface for 
$\operatorname {PSL} (2, {\operatorname {\bold C}})$-representation which
is not necessarily discrete. \par
Assume for a moment that the holonomy representation of the projective
structure is discrete. Then this equivariant map 
actually determines the pleated surface realizing the measured lamination
$\mu$ in the quotient
$3$-manifold for the holonomy representation.  On the other hand, Thurston
gave a remark in
\cite {Th1} that realizing a measured lamination $\mu$ in a hyperbolic
$3$-manifold is a ``harmonic map"  from $[\mu] \in \Cal {PML}$, where
$\Cal {PML}$ is the Thurston boundary of $T_g$. (A rough explanation is
given  in the following way:
$[\mu] \in \Cal {PML}$ is the limit of a degenerating sequence $\{Y_n\}$ of
hyperbolic structures which shrink in the direction $\mu$ as $n \to 
\infty$. Therefore, the harmonic mapping  from $Y_n$ to
a fixed hyperbolic
$3$-manifold stretch along this direction $\mu$ significantly. From the
definition of the energy, the harmonic map from $Y_n$ sends this
direction $\mu$ close to the realization of $\mu$ in the $3$-manifold and
the image is contained in its convex core of the $3$-manifold. Hence 
for large $n$ the image is close to a pleated surface with pleating locus
$\mu$.)
 This intuitive claim was
justified by Minsky
\cite {Mi2}: a pleated surface is the limit (in a very strong sense) of the
images of the harmonic maps from surfaces whose `maximal stretch
direction' is the pleating locus, when the pleating locus is complete.
\par
Since collapsing the grafted part can be seen as a mapping from 
the grafted surface to the pleated surface in the quotient $3$-manifold,
 it is natural to expect collapsing the grafted part is close
to harmonic, when the grafted part is very large.
\par When 
 the holonomy
representation is not discrete, we can still think of harmonic maps 
in the following way:  as in
Donaldson \cite {D}, form a flat
$\operatorname {\bold H}^3$- bundle defined by
$$\Cal {H} = \operatorname {\bold H}^2 \times _\chi \operatorname {\bold
H}^3 \to Y,$$ 
where $\chi$ is the holonomy representation.
Then for a section $s: Y \to \Cal {H}$ take the vertical part of its
derivative: $(Ds)_x : T_x Y \to T_{s(x)}\Cal {H}_x$ where $x$ is a
point on $Y$ and $\Cal {H}_x$ is the fiber over $x$. Then define the
energy by $\Cal {E} (s) = \int _Y ||Ds||^2 dV$ where $dV$ is the volume
form. A {\it twisted harmonic map} is a critical point for the energy
functional. Donaldson \cite {D} showed the existence of the twisted
harmonic map. 
In the same way, we can also define  ``pleated surfaces in
the vertical direction of $\Cal {H}$", which is equivalent to consider
the equivariant map realizing bending. Considering the intuitive
explanation of the relation between harmonic maps and pleated surfaces
for $3$-manifolds, it is reasonable to expect similar things are true
when the representation is not discrete. \par

\endremark
\remark {Remark 2 (an alternative proof)}
When $\mu$ is supported on a simple closed curve,
we can show the properness of the grafting map without using harmonic
maps. In fact, when $\mu$ is supported on a simple closed curve, 
 it is easy to see that in the case (i) the Teichm\"uller
distance between $Y_n$ and $X_n$ is  bounded by a constant independent of
$n$.
 In the case
(ii), we can prove that $ E_{Y_n}(\mu) \ge l_{X_n} (\mu)+ O(1)$ applying
the Thurston metric on $Y_n$ to the definition of the extremal length,
for any $\mu$. However, the author  exhibited the proof using
harmonic maps because it gives a better geometric perspective and also
because it seems (to the author) that for the case (i) arguments by
approximation would not work to give  a uniform constant to bound
 the Teichm\"uller distance between $Y_n$ and $X_n$ for
general measured laminations. 
\endremark
\remark {Remark 3 (properness with respect to $\Cal {ML}$)} Theorem 3.4
implies also that for a fixed hyperbolic surface $X$, the mapping
$gr_{\cdot}(X): \Cal {ML} \to T_g$ is proper.
\endremark
\bigskip
\bigskip

We continue the proof of Theorem 3.2. Although the properness of grafting
map was proved for any measured lamination, we will assume that $\mu$ is
an integral point of $\Cal {ML}$ for the proofs of (2) and (3).
\demo {Proof of (2)}
Here we use the parametrization of projective structures by
$Q_g$, i.e. the space of quadratic differentials.  We will observe how the
fiber of $Q_g$ over each point $Y \in T_g$, i.e. the space of projective
structures on a fixed complex structure,  is mapped by the holonomy map. Let
$Rep = Hom (\pi_1 
\Sigma_g ,
\operatorname {PSL (2,
\bold {C})})/ 
\operatorname {PSL (2, \bold {C})}$ denote the space of 
$\operatorname {PSL (2, \bold {C})}$-representations
of 
$\pi_1 \Sigma_g $, and $hol : Q_g \to Rep$ denote the
holonomy map, namely the mapping which sends each projective structure to
its holonomy representation. \par
 Let $X \in T_g$ be a hyperbolic surface. By the assumption that $\mu$ is 
an integral point, the holonomy
representation of
$Gr_\mu (X)$ is a Fuchsian group
$\Gamma_X$ with quotient surface $X$, hence the holonomy representation  is
in the space of real representations (namely, the equivalence class in $Rep$
with a representative in 
$Hom (\pi_1  \Sigma_g , \operatorname {PSL (2, \bold {R})}$.) which is
denoted by $Rep_{\bold R}$. Now, let  
$Y$ denote the complex structure under the projective structure $Gr_\mu
(X)$ and 
$Q(Y)$ be the space of projective structures on $Y$, i.e. the fiber of
$Q_g$ over $Y \in T_g$. There is an element $\varphi \in Q(Y)$
corresponding to $Gr_\mu (X)$. Then $hol (Q(Y))$ intersects with 
$Rep_{\bold R}$ at $hol (\varphi) = \Gamma_X$. By Faltings' theorem,
(Faltings \cite {F, Theorem 12}), this intersection is transversal.
 Therefore,
at
$hol(Gr_\mu (X)) = \Gamma_X$, we can take a basis $\{u_1, ..., u_{6g-6}\}$
of the (real) tangent space $T_{\Gamma_X}(hol(Q(Y)))$ and a basis $\{v_1,
..., v_{6g-6}\}$ of the (real) tangent space
$T  _{\Gamma_X}(Rep_{\bold R})$ such that 
$\{u_1, ..., u_{6g-6}, v_1, ..., v_{6g-6}\}$ forms the basis of the
tangent space 
$T _{\Gamma_X}(Rep)$. Remember that $hol$ is a local
$C^1$-diffeomorphism (Hejhal). Therefore there are a neighborhood $U$ of
$\varphi$ in $Q_g$ and a neighborhood $V$ of
$\Gamma_X$ in $Rep$ such that $hol|U : U \to V$ is a $C^1$-diffeomorphism
and the inverse map $g:V \to U$ of $hol|U$ is well-defined. By the
bundle structure of $Q_g$, we can take a neighborhood $U'$ of $Y$ in $T_g$ 
such that the restriction $\pi|\pi^{-1}(U') : \pi^{-1}(U') \to U'$
is identified with the product of $U'$ with  $\operatorname {\bold
R}^{6g-6}$, where $\pi : Q_g \to T_g$ is the projection. Then we may
assume that 
$U$ is the product of $U'$ and  an open set of
$\operatorname {\bold R}^{6g-6}$. Denote by $\phi$ the point of
$\operatorname {\bold R}^{6g-6}$ such that $\varphi$ corresponds to 
$(Y, \phi) \in T_g \times \operatorname {\bold R}^{6g-6}$.
Then the tangent space
$T _\varphi U $ is spanned by the `direction of the base space'
$T _Y T_g $ and the `direction of the fiber'
$T_\phi \operatorname {\bold R}^{6g-6}$ and the derivative
$dg$ maps $T_{\Gamma} Rep$ onto $T_\varphi T_g$. Now, $\{dg (u_i)\}_{i =
1,...,6g-6}$ is contained in the direction of fiber. Therefore, none of
the non-zero vectors in $dg (T_{\Gamma_X}Rep_{\bold R})$
is contained in the direction of fiber. It follows that $d(\pi \circ g) :
T_{\Gamma_X} Rep_{\bold R}  \to  T _Y T_g $
is surjective. As we can identify the component of
the space of real representations containing $\Gamma_X$ with the
Teichm\"uller space, the composition of the restriction
$g|Rep_{\operatorname {\bold R}}$ with $\pi$ is equal to $gr_\mu$. Therefore,
$gr_\mu$ is locally diffeomorphic at $X$.
\qed
\enddemo

\demo {Proof of (3)} 

Let $QF(\mu)$ be the set of projective structures obtained by  
quasiconformal deformations of a grafted projective structure $Gr_\mu
(X)$.  Then
$QF(\mu)$ is identified with  the space of quasiconformal deformations of 
the holonomy representation of
$Gr_\mu (X)$, which is a Fuchsian group (cf. \cite {ST}). Recall that
$T_g$ has a natural complex structure and the space of quasiconformal
deformations of the Fuchsian group is identified with the complex
manifold $T_g \times T_g$. With respect to this identification, the
mapping
$\Pi : QF(\mu) \to T_g$ which sends each projective structure in
$QF(\mu)$ to the underlying complex structure is holomorphic (cf. \cite
{ST}). Now, in the space of quasiconformal deformations of a Fuchsian
group, the set of Fuchsian groups , which is identified with $T_g$,
forms a real analytic submanifold. Therefore, the restriction of $\Pi :
QF(\mu)
\to T_g$ to this set of Fuchsian groups is real analytic. This
restriction is the same mapping as $gr_\mu :T_g \to T_g$.
\qed
\enddemo

\enddocument